# Smooth Global Solution of the Navier-Stokes Equation

Maoting Tong，Public Basic Department of Nanjing Vocational University of Industry Technology, Nanjing, 210023, P. R. China. (e-mail: 3977093889@qq.com).

Daorong To$n$, Corresponding author, Department of Mathematics of Hohai University, Nanjing, 210098, P. R. China. (e-mail: 1760724097@qq.com).

## ABSTRACT

Whether the Navier-Stokes equation has a Solution depends on its initial state. Suppose that the given externally applied force $f(t,x)=0$. In the framework function space $L_p$ $(R^3)$, the Navier-Stokes equation has a unique smooth global solution with bounded energy if and only if the initial value $u_0$ is smooth and divergence free vector field.

## I. INTRODUCTION

Suppose the given externally applied force $f(t,x)=0$， the Navier-Stokes initial value problem

can be written in the following form

$$\begin{cases} \frac{\partial u}{\partial t} = \Delta u - \nabla p - (u \bullet \nabla)u, x \in R^3, t \in [t_0, \infty), \\ \nabla \bullet u \equiv \operatorname{div} u = 0, \\ u|_{t=t_0} = u_0\,, x \in R^3 \end{cases} \tag{1}$$

where $u=u(t,x)=(u_1(t,x),u_2(t,x),u_3(t,x))$ is the velocity field, $u_0=u_0(x)$ is the initial velocity,

$p=p(t,x)$ is the pressure, and

$$\Delta u = (\sum_{i=1}^{3}\frac{\partial^2 u_1}{\partial x_i^2}, \sum_{i=1}^{3}\frac{\partial^2 u_2}{\partial x_i^2}, \sum_{i=1}^{3}\frac{\partial^2 u_3}{\partial x_i^2}).$$

The Navier-Stokes equation describes the momentum conservation of viscous incompressible fluid. This equation generalizes the general law of a viscous incompressible fluid flow. Its solution does not satisfy the linear superposition principle, because viscous incompressible fluid flow is a nonlinear system. The equation of motion of a viscous fluid was first proposed by Navier in 1827 and only considers the flow of an incompressible fluid. Poisson proposed an equation of motion for a compressible fluid in 1831. In 1845, Saint Venant and Stokes independently proposed that the viscosity coefficient is a constant form, which is called the Navier-Stokes equation, or the N-S equation for short. Its solution is challenging and complex. Before further development and breakthrough of the solution idea or technology, the accurate solution could only be obtained for some simple special flows. However, in some cases, the equation can be simplified to obtain an approximate solution. For example, when the Reynolds number $Re > 1$, the viscous force outside the boundary layer of the flowing object is far less than the inertial force, the viscous term in the equation can be ignored, and the N-S equation can be simplified to the Euler equation in the ideal flow. In the boundary layer, the N-S equation can be simplified to a boundary layer equation.

The existence, uniqueness and regularity of solutions to N-S equation have been extensively studied. There is much literature on the solvability of the initial value problem for N-S equations. The earliest research result on this problem can be found in a series of remarkable papers [13] [14] by J. Leray, which appeared approximately 90 years ago and are the pioneering work about N-S equation. In the 1960s Japanese mathematicians H. Fujita and T. Kato deduced the existence theorem for three-dimensional flows through a Hilbert space approach, employing theories of fractional power and semigroups of operators. This work is a groundbreaking work. Philip Isett's [9] and T.Kato's [10] are all good jobs, too.

The existence of smooth solutions of N-S equations in three-dimensional space was set as one of the seven Millennium Prize Problems by the Clay Institute of Mathematics in the United States. C. L. Fefferman put forward the following specific requirements for this problem.

Existence and smoothness of Navier-Stokes solutions on $R^3$:

Let $u_0$ be any smooth, divergence free vector field satisfying

$$|\partial_x^\alpha u_0(x)| \le C_{\alpha k}(1+|x|)^{-k} \text{ on } R^3 \qquad (*)$$

for any $\alpha$ and $k$. Take the given externally applied force $f(t,x)$ to be identically zero. Then there exist smooth functions $p(t,x), u(t,x)$ on $[0,\infty)\times R^3$ that satisfy (1) and

$$p, u \in C^\infty([0,\infty)\times R^3),$$

$$\int_{R^3} |u(t,x)|^2\, dx < C \text{ for all } t \ge 0 \text{(bounded energies).}$$

There are different views of the N-S equation in the international mathematical community. In 2015, Terence Tao constructed an equation nearly equivalent to the original N-S equation in an 76-page article《Finite time blowup for an averaged three-dimensional Navier-Stokes equation》(arXiv:1402. 0290）, and proved that the equivalent equation would explodes in a finite time. It can then be inferred that the original N-S equation will also explode in finite time. In other words, he does not believe that the N-S equation will have a smooth solution. In 2017, Tristan Buckmaster (Princeton) and Vlad Vicol (NYU) in the paper《Non-uniqueness of weak solutions to the Navier-Stokes equation》（arXiv:1709.10033）proved the non-uniqueness of weak solutions to the Navier-Stokes equation, thus questioning its ability to fully express a hydrodynamic system. However, they did not substantively prove that the N-S equation explodes in finite time or that it has multiple solutions. Today, some mathematicians tend to believe that the N-S equation has a unique smooth solution. The focus of the problem is what types of solutions will be available.

It is difficult to solve the third Millennium Problem. Fefferman said: "Standard methods from PDE appear inadequate to settle the problem. Instead, we probably need some deep, new ideas" in the last sentence of his paper. The theory of semigroups of operators produces new results and involves. In the theory of semigroups of operators we have involved two types of structures in the three major structures of mathematics, topological structure (convergence) and algebraic structure (semigroups of operators). But this is not enough to solve a difficult mathematical problem in history. We also need the third kind of important mathematical structure—order structure, that is, the theory of Banach lattice. Thus, we could reveal all the secrets of N-S equations. We provide a sufficient and necessary condition for the N-S equation having a unique smooth global solution. We prove that the N-S equation has a smooth global solution with bounded energy if and only if the initial value $u_0$ is smooth, divergence free vector field. We followed the methods of Maoting Tong at [22] [23][24].

In this study the range of integers $p$ is $2 \le p<\infty$. Let $L_p(R^3)$ be the Banach space of real vector functions with the inner product defined in the usual manner, that is,

$$L_p(R^3) = \{u: R^3 \to R^3, u = (u_1, u_2, u_3), u_i \in L^p(R^3)(i = 1,2,3)\}.$$

Here, superscript $p$ represents the space of the scalar-valued value functions, and subscript $p$ represents the space of the vector value functions.

All the functions we discuss in this article belong to $L_p(R^3)$. In the framework function space $L_p(R^3)$, the set of all real vector functions $u$ such that $\operatorname{div} u = 0$ and $u \in C_0^\infty(R^3)$ is denoted by $C_{0,\sigma}^\infty(R^3)$.

**Definition 1.** Let $DL_p(R^3)$ denote the closure of $C_{0,\sigma}^\infty(R^3)$ in $L_p(R^3)$. This function space is called quasi-divergence free space.

We choose the quasi-divergence free space as the function space to solve the N-S equations. The reason for this choice is that taking the given externally applied force f (t, x) to be identically zero. In this way, the only internal force is pressure, so the velocity field does not generate vortices, and the velocity field is smooth and nearly irrotational. The function space of solutions should be a complete Banach space, and the quasi-divergence free space just meets these requirements.

If $u \in C^\infty(R^3)$ then $u \in DL_2(R^3)$ implies $div\, u = 0$ (see p.270 and p.313 in [5]). Now we prove that if $u \in C^\infty(R^3)$ then

$$u \in DL_p(R^3) \text{ implies } \operatorname{div} u = 0. \tag{2}$$

In fact, suppose that $u = (u_1, u_2, u_3) \in C^\infty(R^3) \cap DL_p(R^3)$. Then for any $x \in \mathrm{R}^3$, there exists a closed set $\Omega \subset R^3$ such that $x \in \Omega$. Then $u = (u_1, u_2, u_3) \in C^\infty(\Omega) \cap DL_p(\Omega)$. Since $DL_p(\Omega) = \overline{C_{0,\sigma}^\infty}$ in $L_p(\Omega)$, $DL_2(\Omega) = \overline{C_{0,\sigma}^\infty}$ in $L_2(\Omega)$ *and* $L_p(\Omega) \subset L_2(\Omega)$, hence $DL_p(\Omega) \subset DL_2(\Omega)$. $u \in DL_p(\Omega)$ implies $u \in DL_2(\Omega)$. From p.270 in [5] we observe that $\operatorname{div} u = 0$ for $x \in \Omega$. Therefore $\operatorname{div} u = 0$ for all $x \in R^3$.

Regarding the relationship between the above spaces, we have

$$C_{0,\sigma}^\infty(R^3) \subset C_0^\infty(R^3) \subset L_p(R^3),$$

$$DL_p(R^3) = \overline{C_{0,\sigma}^\infty(R^3)} \subset L_p(R^3),$$

and

$$DL_p(R^3) \subset L_p(R^3) = W_{0,p}(R^3),$$

$$\|\bullet\|_{DL_p(R^3)} = \|\bullet\|_{L_p(R^3)}.$$

From Theorem 3(a) in [7, p.129] we have the Helmholtz decomposition

$$L_p(R^3) = DL_p(R^3) \oplus DL_p(R^3)^\perp$$

where $DL_p(R^3)^{\perp} = \{\nabla h; h \in W^{1,p}(R^3)\}$. Let $P$ be the orthogonal projection of $L_p(R^3)$ onto $DL_p(R^3)$. Applying $P$ to the first equation of (1) and considering the other equations, we obtain an abstract initial value problem,

$$\begin{cases} \frac{du}{dt} = P\Delta u + Fu, t \in (0,\infty), \\ u|_{t=0} = u_0, \ x \in R^3, \end{cases} \qquad (3)$$

where $Fu = -P(u \bullet \nabla)u$.

The integral form of (3) is

$$u(t) = e^{tP\Delta}u_0 + \int_{t_0}^{t} e^{(t-s)P\Delta} Fu(s)ds. \qquad (4)$$

Fujita and Kato discussed the N-S initial value problem in the finite case $[0,T] \times \Omega$ where $\Omega$ is a closed domain in $R^3$. In Theorem 1.6 of [5] they proved that if $u$ is a solution to (3), then $u$ is of class $C^{\infty}$as $L_2(\Omega)$ −valued function. They provided a sufficient condition for (3) to obtain a solution. Giga and Miyakawa discussed the solutions to the N-S initial value problem in $L_r(1<r<\infty)$. In Theorem 3.4 of [6] they proved that the solution to (4) belongs to $(C^{\infty}(\overline{\Omega} \times (0,T]))^n$. They provided a sufficient condition for (3) and (4) to obtain a solution.

## II. SOME LEMMAS

We first give some lemmas. Lemmas 1 and 2 prove that the Laplace operator $\Delta$ limited to $DL_p(R^3)$ is the infinitesimal generator of an analytic semigroup of contractions on $DL_p(R^3)$. In this way, we can use the semigroups tool to discuss the N-S initial value problem. Lemma 3 gives the estimation inequality of the nonlinear term $(u \bullet \nabla)u$ to deal with the nonlinear term in the N-S equation. Lemma 4 gives a sufficient condition for the global solution of the initial value problem (3), which is a bridge for deriving the main results. Lemma 5 uses the Banach lattices to prove the closeness of the binary operation $(u \bullet \nabla)v$ in $DL_p(R^3)$.

Let $u = (u_1, u_2, u_3) \in L_p(R^3)$,the operator

$$-\Delta u_i = -\sum_{j=1}^{3} \frac{\partial^2 u_i}{\partial x_j^2} (u_i \in L^p(\Omega), i = 1,2,3)$$

is strongly elliptic of order 2 in $L^p(\Omega)$. From formula (2.3) in [1] and a more general result in [2, p.176], Theorem 7.3.2 in [18] holds for $\Omega = R^3$. Using a proof method completely similar to Theorem 7.3.6 in [18], we can prove that $\Delta$ is the infinitesimal generator of an analytic semigroup of contractions on $L^p(R^3)$ with $D(\Delta) = W^{2,p}(R^3) \cap W_0^{1,p}(R^3)$ (see formula 7.3.5 in [18]). Hence, $\Delta$ is also the infinitesimal generator of an analytic semigroup of contraction on $L_p(R^3)$ with $D(\Delta) = W_{2,p}(R^3) \cap W_{1,p,0}(R^3)$, where $W_{2,p}(R^3)$ and $W_{1,p,0}(R^3)$ are the Sobolev spaces of vector value in $W^{2,p}(R^3)$ and $W_0^{1,P}(R^3)$ respectively.

**Lemma 1.** For every $u \in L_p(R^3)$, $div\, u = 0$ if and only if $div\,[(\lambda I - \Delta)]u = 0$ for

$\lambda \in \Sigma_{\vartheta} = \{\lambda: \vartheta - \pi < \arg\lambda < \pi - \vartheta, |\lambda| \geq r\}$

where $0<\vartheta< \frac{1}{2}$, and $r$ is a positive number.

Proof. Let $u = (u_1, u_2, u_3) \in L_p(R^3)$. Then

$$\Delta u=\left(\frac{\partial^2u_1}{\partial {x_1}^2}+\frac{\partial^2u_1}{\partial {x_2}^2}+\frac{\partial^2u_1}{\partial {x_3}^2},\frac{\partial^2u_2}{\partial {x_1}^2}+\frac{\partial^2u_2}{\partial {x_2}^2}+\frac{\partial^2u_2}{\partial {x_3}^2},\frac{\partial^2u_3}{\partial {x_1}^2}+\frac{\partial^2u_3}{\partial {x_2}^2}+\frac{\partial^2u_3}{\partial {x_3}^2}\right),$$

$$\mathrm{div}(\Delta u)=\frac{\partial^3u_1}{\partial {x_1}^3}+\frac{\partial^3u_1}{\partial x_1\partial {x_2}^2}+\frac{\partial^3u_1}{\partial x_1\partial {x_3}^2}+\frac{\partial^3u_2}{\partial x_2\partial {x_1}^2}+\frac{\partial^3u_2}{\partial {x_2}^3}$$

$$+\frac{\partial^3u_2}{\partial x_2\partial {x_3}^2}+\frac{\partial^3u_3}{\partial x_3\partial {x_1}^2}+\frac{\partial^3u_3}{\partial x_3\partial {x_2}^2}+\frac{\partial^3u_3}{\partial {x_3}^3}=\frac{\partial^2}{\partial {x_1}^2}(\frac{\partial u_1}{\partial x_1}+\frac{\partial u_2}{\partial x_2}+\frac{\partial u_3}{\partial x_3})+\frac{\partial^2}{\partial {x_2}^2}(\frac{\partial u_1}{\partial x_1}+\frac{\partial u_2}{\partial x_2}$$

$$+\frac{\partial u_3}{\partial x_3})+\frac{\partial^2}{\partial {x_3}^2}(\frac{\partial u_1}{\partial x_1}+\frac{\partial u_2}{\partial x_2}+\frac{\partial u_3}{\partial x_3})$$

$$=\left(\frac{\partial^2}{\partial {x_1}^2}+\frac{\partial^2}{\partial {x_2}^2}+\frac{\partial^2}{\partial {x_3}^2}\right)\left(\frac{\partial u_1}{\partial x_1}+\frac{\partial u_2}{\partial x_2}+\frac{\partial u_3}{\partial x_3}\right)$$

$$=\Delta(\mathrm{div}\,u).$$

Therefore, we have

$$\mathrm{div}[(\lambda I-\Delta)u]=(\lambda I-\Delta)(\mathrm{div}u), \qquad (5)$$

from which it is clear that $div\ u=0$ implies that $\mathrm{div}\,[(\lambda I-\Delta)u]=0$.

Since the Theorem 7.3.2 in [18] holds for $\Omega=R^3$, there exist constant $C>0$, $r\geq 0$ and $0<\vartheta<{}^{\pi}/_{2}$ such that

$$\|u\|_{L_p(R^3)}\leq\frac{C}{|\lambda|}\|(\lambda I-\Delta)u\|_{L_p(R^3)} \qquad (6)$$

for $u\in W_{2,p}(R^3)\cap W_{1,p,0}(R^3)$ and $\lambda\in\Sigma_\vartheta=\{\lambda:\vartheta-\pi<\arg\lambda<\pi-\vartheta,|\lambda|\geq r\}$. From (6), it follows that for every $\lambda\in\Sigma_\vartheta$, $\lambda I-\Delta$ is injective from $D(\Delta)$ into $L_p(R^3)$. From (5), it follows that $\mathrm{div}[(\lambda I-\Delta)u]=0$ implies that $div\ u=0$.☐

In [18], the proof of Theorem 1.7.7 used the Theorem 1.5.2, whose condition is only $\rho(A)\supset R^+=\{t\in R:t>0\}$. $\rho(A)$ *is the resolvent set of* $A$. Therefore, the condition $(i)$ of Theorem 1.7.7 can be changed to $\rho(A)\supset\Sigma_\delta=\{\lambda:|arg\ \lambda|<\pi/2+\delta\}$ *for some* $0<\delta<\pi/2$. The proof of Theorem 2.5.2 in [18] used the Theorem 1.7.7. Therefore, $\cup\{0\}$ can be deleted in the formula (2.5.4) of Theorem 2.5.2 in [18] and the formula can be changed to

$$\rho(A)\supset\textstyle\sum=\{\lambda:|arg\ \lambda|<\pi/2+\delta\}$$

for some $0<\delta<\pi/2$. In this way, we can refer to Lemma 2.5.2 from [18] in the proof of following Lemma 2.

**Lemma 2.** The operator $\Delta_{|DL_p(R^3)}$ is the infinitesimal generator of an analytic semigroup of contractions on $DL_p(R^3)$.

Proof. The proof of Theorem 7.3.5 in [18] used formula (7.3.4) in [18], which as pointed out earlier, holds for $\Omega=R^3$, and hence Theorem 7.3.5 in [18] is also valid for $\Omega=R^3$. From this theorem, $\Delta$ is an infinitesimal generator of an analytic semigroup on $L^p(R^3)$. Therefore, $\Delta$ is also an infinitesimal generator of an analytic semigroup on $L_p(R^3)$. Let $\{T(t)|t\geq 0\}$ be the restriction of the analytic semigroup generated by $\Delta$ on $L_p(R^3)$ to the real axis. Thus, $\{T(t)|t\geq 0\}$ is a $C_0$ −semigroup of contractions. We have already noted that $DL_p(R^3)$ is a closed subspace of $L_p(R^3)$.

For every $u \in L_p(R^3)$ with $div\, u = 0$ and $\lambda \in \rho(\Delta) \cap \Sigma_\vartheta = \{\lambda: \vartheta - \pi < \arg\lambda < \pi - \vartheta, |\lambda| \geq r\}$ we have $(\lambda I - \Delta)[R(\lambda:\Delta)u] = u$, where $\sum_\vartheta$ is the same as in the proof of Lemma 1. From Lemma 1, it follows that $divR(\lambda:\Delta)u = 0$, i.e., $R(\lambda:\Delta) = (\lambda I - \Delta)^{-1}$ is preserving the divergence free property for $\lambda \in \rho(\Delta) \cap \Sigma_\vartheta$. From Theorem 2.5.2 (c) in [18] it follows that $\rho(\Delta) \supset R^+$, and so $\rho(\Delta) \cap \Sigma_\vartheta \supset \{\lambda: \lambda \geq r\}$. Hence $R(\lambda:\Delta)$ is preserving the divergence-free property for every $\lambda \geq r$. Let $u \in DL_p(R^3)$. Then there exists a sequence $u_n$ such that $\lim_{n\to\infty} u_n = u$ and div $u_n = 0$ for $n = 1,2,\ldots$. Since $R(\lambda:\Delta)$ is bounded and so is continuous. Hence, $\lim_{n\to\infty} R(\lambda:\Delta)u_n = R(\lambda:\Delta)u$and $divR(\lambda:\Delta)u_n = 0$ for every $\lambda \geq r$. Therefore, $R(\lambda:\Delta)u \in DL_p(R^3)$ for every $\lambda \geq r$. It follows that $DL_p(R^3)$ is $R(\lambda:\Delta)$-invariant for every $\lambda \geq r$. Now, Theorem 4.5.1 in [18] implies that $DL_p(R^3)$ is $T(t)_{t\geq 0}$ $-$invariant. From 1.5.12 in [11] and Proposition 2.2.3 in [11], it follows that $\Delta_{|DL_p(R^3)}$ is the infinitesimal generator of the $C_0$- semigroup $\left\{T(t)_{|DL_p(\Omega)}: t \geq 0\right\}$ of contractions on $DL_p(R^3)$.

Suppose that $\lambda \in \rho(\Delta)$, i.e., there exists $R(\lambda:\Delta)$ from $L_p(R^3)$ into $D(\Delta)$. Then, for any $u \in DL_p(R^3) \subset L_p(R^3)$, $R(\lambda:\Delta)u \in DL_p(R^3)$ because $R(\lambda:\Delta)$ preserves divergence-free property,

we have

$$(\lambda I - \Delta)R(\lambda:\Delta)u = u,$$

$$R(\lambda:\Delta)(\lambda I - \Delta)u = u. \tag{7}$$

Thus, the formula (7) becomes $(\lambda I - \Delta_{|DL_p(R^3)})R(\lambda:\Delta)_{|DL_p(R^3)}u = u$ and

$$R(\lambda:\Delta)_{|DL_p(R^3)}(\lambda I - \Delta_{|DL_p(R^3)})u = u.$$

Hence, $(\lambda I - \Delta_{|DL_p(R^3)})R(\lambda:\Delta)_{|DL_p(R^3)} = I$ and $R(\lambda:\Delta)_{|DL_p(R^3)}(\lambda I - \Delta\big|DL_p(R^3)) = I$. We obtain

$$(\lambda I - \Delta_{|DL_p(R^3)})^{-1} = R(\lambda:\Delta)_{|DL_p(R^3)}. \tag{8}$$

From formula (8) and Theorem 2.5.2(c) in [18] we have

$$\rho(\Delta_{|DL_p(R^3)}) \supset \rho(\Delta) \supset \Sigma = \left\{\lambda: |arg\, \lambda| \langle {}^{\pi}/_{2} + \delta\right\},$$

where $0<\delta<{}^{\pi}/_{2}$. Thus, for $\lambda \in \Sigma$, $\lambda I - \Delta_{|DL_p(R^3)}$is invertible. From Theorem 2.5.2 (c) in [18] we have, for $\lambda \in \Sigma, \lambda \neq 0$,

$$\left\|R(\lambda:\Delta_{|DL_p(R^3)})\right\| = \underset{\substack{u\in DL_p(R^3)\\ \|u\|_{DL_p(R^3)}=1}}{Sup} \left\|R(\lambda:\Delta_{|DL_p(R^3)}|)u\right\|$$

$$\leq \underset{\substack{u\in L_p(R^3)\\ \|u\|_{L_p(R^3)}=1}}{Sup} \|R(\lambda:\Delta)u\| = \|R(\lambda:\Delta)\|_{L_p(R^3)} \leq \frac{M}{|\lambda|}. \tag{9}$$

Theorem 2.5.2(c) in [18] implies that $\left\{T(t)_{|DL_p(R^3)}: t \geq 0\right\}$ can also be extended to analytic semigroup on $DL_p(R^3)$. Therefore, $\Delta_{|DL_p(R^3)}$ is an infinitesimal generator of an analytic semigroup of contractions on $DL_p(R^3)$. □

Suppose that $-A$ is the infinitesimal generator of an analytic semi-group $T(t)$ on a Banach space $X$. We can define the fraction powers $A^{\alpha}$ for $0 \le \alpha \le 1$, and $A^{\alpha}$ is a closed linear invertible operator with domain $D(A^{\alpha})$ dense in $X$. $D(A^{\alpha})$, equipped with the norm $|||u||| = \|u\| + \|A^{\alpha}u\| (u \in X)$ which is equavalent to $\|A^{\alpha}u\|$, is a Banach space denoted by $X_{\alpha}$. For example, if $X = DL_p(R^3)$ and $-A = \Delta$, then the Banach space $DL_p(R^3)_{1/2}$ is $D((-\Delta)^{\frac{1}{2}})$ equipped with the norm $\|u\|_{1/2} = \left\|(-\Delta)^{\frac{1}{2}}u\right\|$.

**Definition 2.** $DL_p(R^3)_{1/2}$ is called the graph space of the fractional power of order 1/2 of the strong elliptic operator $-\Delta$ in $DL_p(R^3)$.

The inclusion relationships between the above function spaces are shown in the figure 1 below.

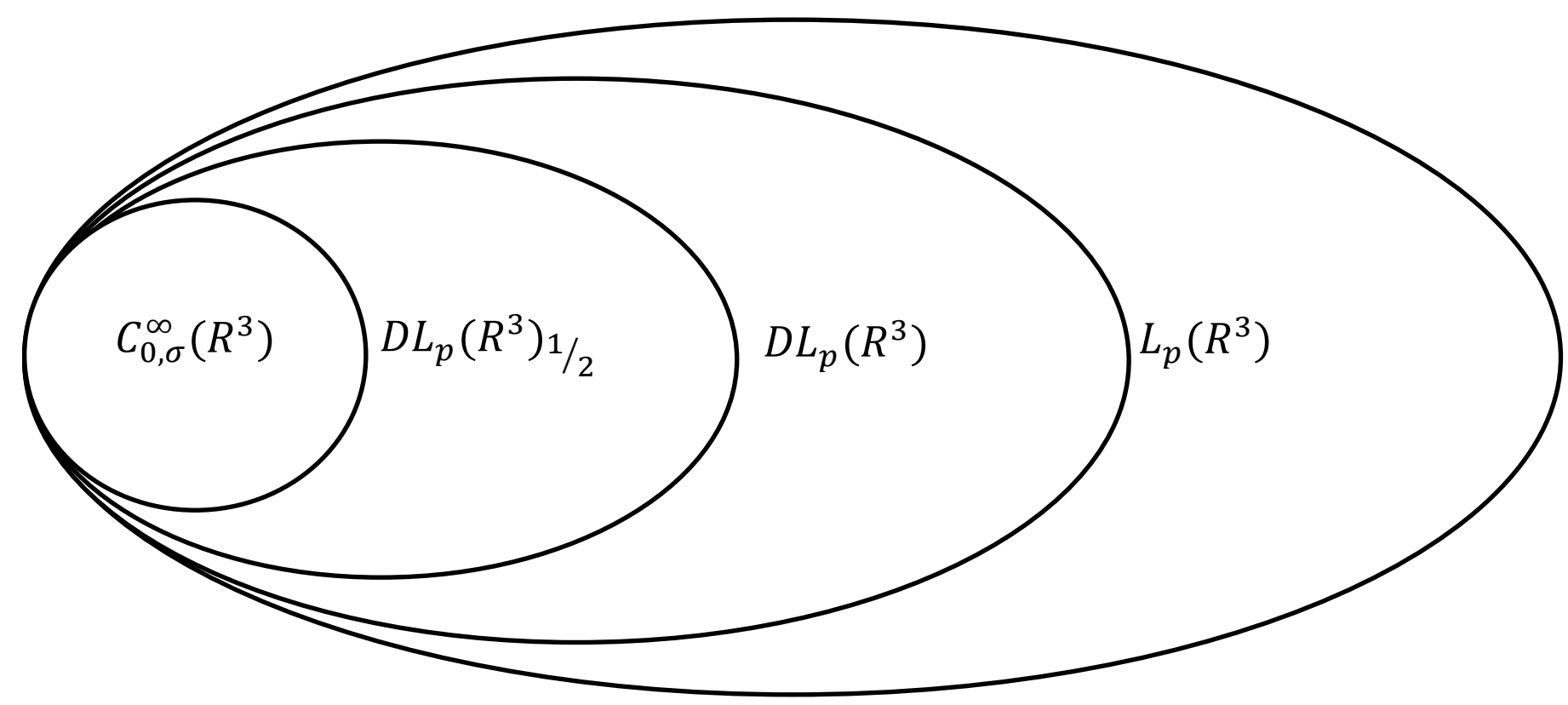

figure 1

In [5] let $C_{0,\sigma}^{1}(\Omega)$ denote the set of all solenoidal functions in $C_0^1(\Omega)$ *where* $\Omega$ *is a bounded* domain in $R^3$. *Let* $H_{0,\sigma}^{1}(R^3)$ be the completion of $C_{0,\sigma}^{1}(\Omega)$ in $C_0^1(\Omega)$. Then $D((-P\Delta)^{1/2}) = H_{0,\sigma}^{1}(\Omega)$ where $P$ is the orthogonal projection from $L_2(\Omega)$ onto $DL_2(\Omega)$ (See $p. 270_4$ in [5]) and

$$C_{0,\sigma}^{\infty}(\Omega) \subset c_{0,\sigma}^{1}(\Omega) \cap DL_2(\Omega) \subset H_{0,\sigma}^{1}(\Omega) \cap DL_2(\Omega)$$

$$D((-P\Delta)^{1/2}) \cap DL_2(\Omega) = D((-\Delta)^{1/2}) \cap DL_2(\Omega)$$

because when $u \in DL_2(\Omega)$, then $\Delta u \in DL_2(\Omega)$ (Lemma 2, it can also be directly proven) and so $P\Delta u = \Delta u$.

In our frame $DL_p(R^3)$, we have similarly

$$C_{0,\sigma}^{\infty}(R^3) \subset C_{0,\sigma}^{1}(R^3) \cap DL_p(R^3)$$

$$\subset D((-\Delta)^{1/2}) \cap DL_p(R^3) = DL_p(R^3)_{1/2}. \tag{10}$$

For $u \in DL_p(R^3)$ we have

$$\|\nabla u\|_{DL_p(R^3)} = \|-\nabla u\|_{DL_p(R^3)} = \left\|(-P\Delta)^{1/2}u\right\|_{DL_p(R^3)}$$

$$= \left\|(-\Delta)^{1/2}u\right\|_{DL_p(R^3)} = \|u\|_{DL_p(R^3)_{1/2}}. \tag{11}$$

(See the similar formula (1.5)′ in [5]. For $u \in DL_p(R^3), -P\Delta u = -\Delta u$.)

Giga proved the following result (Lemma 2.2 in [5]): Let $0 \le \delta < {}^{1}/_{2}\,|-n(1-r^{-1})/2|, 1<r<\infty, n \in N$. Then

$$\left\|A^{-\delta}P(u,\nabla)v\right\|_{0,r} \le M\left\|A^{\vartheta}u\right\|_{0,r}\|A^{\omega}v\|_{0,r}$$

with some constant $M = M(\delta,\vartheta,\omega,r)$, provided $\delta+\vartheta+\omega \ge n/2r+1/2, \vartheta>0, \omega>0$, $\omega+\delta \ge 1/2$. (The original text requires $\omega+\sigma>1/2$ , but a careful look at the proof shows that this condition is not used at all, can be changed to $\omega+\sigma \ge 1/2$.)

From this result , if we take $n=1, r=p$, $\delta=0$, and $\theta=\omega=1/2$, then we have

$$\|(u\bullet\nabla)v\|_{DL_p(R^3)} = \|(u\bullet\nabla)v\|_{L_p(R^3)}$$

$$= \|(u\bullet\nabla)v\|_{W^{0,p}(R^3)}$$

$$\le M\left\|(-\Delta)^{\vartheta}u\right\|_{DL_p(R^3)}\|(-\Delta)^{\omega}v\|_{DL_p(R^3)}$$

$$\le M\|u\|_{DL_p(R^3)_{1/2}}\|v\|_{DL_p(R^3)_{1/2}}$$

with some constants $M$ for any $u,v \in DL_p(R^3)$.

In the following Lemma 4 we will prove that $u,v \in DL_p(R^3)$imply $(u\cdot\nabla)v \in DL_p(R^3)$.Hence $P(u\cdot\nabla)v = (u\cdot\nabla)v$ and we have

**Lemma 3.** Suppose that $u$, $v \in DL_p(R^3)$, then

$$\|(u\bullet\nabla)v\|_{DL_p(R^3)} \le M\,\|u\|_{DL_p(R^3)_{1/2}}\|v\|_{DL_p(R^3)_{1/2}}.$$

A. Pazy said in his book [18, p.69] : " Most of the results on fractional powers that we will obtain in this section remain true if $0 \notin \rho(A)$". In [18], the proof of Theorem 6.3.1 used Theorem 2.6.13, the proof of Theorem 2.6.13 used the formula (2.6.7), the formula (2.6.7) follows from the formula (2.6.6), the formula (2.6.6) follows from the formula (2.5.20), the formula (2.5.20) is equivalent to the formula (2.5.6) of Theorem (2.5.2), the proof of the formula (2.5.6) used Theorem (1.7.7), the proof of Theorem (1.7.7) used Theorem (1.5.2), the proof of Theorem (1.5.2) used only the condition that the resolvent set contains positive real number. That is, the condition $0 \in \rho(A)$ can be deleted in Theorem 1.7.7, the conditions $0 \in \rho(A)$ can be deleted in Theorem (2.5.2), the condition $0 \in \rho(A)$ can be deleted in Theorem (2.6.13), and so the condition $0 \in \rho(A)$ can be deleted in Theorem 6.3.1. Likely, the proof of following Lemma 4 used also Theorem 2.6.13. Therefore, there is no condition $0 \in \rho(A)$ in the following Lemma 4. Lemma 4 extends Theorem 6.3.1 in [18] to the case in which both time and space are infinite, and the proof approach is similar.

Although $0 \notin \rho(\Delta)$, we can use the formula (2.6.9) in [18] to define the fractional power of Δ. Y. Giga and T. Miyakawa used $(-\Delta)^{-1/2}$ in [6].

The proof of Theorem 6.3.1 in [18] does not use the condition that U is an open set, which can be omitted. In the following assumption, we only assume that U is a subset, not necessarily is an open set.

**Assumption (F).** Let $X = DL_p(R^3)$, and let $U$ be a subset of $R^+ \times X_\alpha$ $(0 \le \alpha \le 1)$. The function $f: U \to X$ satisfies the assumption (F) if for every $(t,u) \in U$ there is a subset $V \subset U$and constants $L \ge 0, 0<\vartheta \le 1$ such that $(t,u) \in V$ and for all $(t_i,u_i) \in V$ $(i = 1,2)$,

$$\|f(t_1,u_1) - f(t_2,u_{2)}\| \le L(|t_1 - t_2|^{\vartheta} + \|u_1 - u_2\|_{\alpha}) \tag{12}$$

**Lemma 4.** Let $-A$ be the infinitesimal generator of an analytic semigroup $T(t)$ on the Banach space $X = DL_p(R^3)$ satisfying $\|T(t)\| \le 1$. If, $\alpha = \frac{1}{2}$ and $f$ satisfies the assumption $(F)$, then for any initial data $(t_0,u_0) \in U$, the equation

$$\begin{cases} \frac{du(t)}{dt} + Au(t) = f(t,u(t)), t \in (t_0,\infty) \\ \qquad\qquad u(t_0) = u_0 \end{cases} \tag{13}$$

with $f(t,u(t)) = -(u \bullet \nabla)u$ has a global solution

$$u \in C([t_0,\infty): DL_p(R^3)) \cap C^1((t_0,\infty): DL_p(R^3)).$$

Proof. (Step 1) From the assumption of operator $A$, it follows that

$$\|A^{\alpha}T(t)\| \le C_{\alpha}t^{-\alpha} \tag{14}$$

for $t \in (t_0,\infty)$ by Theorem 2.6.13(c) in [18]. For fixed $(t_0,u_0) \in U$, we choose $\delta > 0$ such that estimate (12) with some constant $L$ and $\vartheta$ holds in the set $V = \{(t,u): t_0 \le t,\ \|u - u_0\| \le \delta\}$. Choose $t_1' > 0$, and let

$$B = \max_{0 \le t \le t_1'} \|f(t,u_0)\|_X. \tag{15}$$

Since $\lim_{t\to 0}\|T(t)\|_{B(DL_p(R^3))} = \|I\|$ implies $\lim_{t\to 0}\left\|\frac{1}{t+1}T(t)\right\|_{B(DL_p(R^3))} = \|I\|$ where $B(DL_p(R^3))$ is the Banach space of all bounded linear operators from $DL_p(R^3)$ into $DL_p(R^3)$, we can choose $t_1$ such that

$$\left\|\frac{1}{t+1}T(t)A^{\alpha}u_0 - A^{\alpha}u_0\right\|_{DL_p(R^3)} < 1/2 \text{ for } t_0 \le t < t_1 \tag{16}$$

and

$$0 < t_1 - t_0 < \min\left\{t_1' - t_0, \left[\frac{\delta}{2}(1-\alpha)C_{\alpha}^{-1}(B + \ \delta L)^{-1}\right]^{\frac{1}{1-\alpha}}\right\}. \tag{17}$$

(Step 2) Let $Y$ be the space of functions $y(t)$ from $[t_0,\infty)$ to $DL_p(R^3)$ which satisfies the following conditions:

1. $y(t)$ is continuous for $t \in [t_0,\infty)$, i.e., $\lim_{t\to t_1}\|y(t) - y(t_1)\|_{DL_p(R^3)} = 0$ for any $t_1 \in [t_0,\infty)$;

2. $\int_{t_0}^{\infty}\|y(t)\|_{DL_p(R^3)} dt$ converges, so $\|y(t)\|_{DL_p(R^3)}$ is bounded for $t \in [t_0,\infty)$.

(1) For $y(t) \in Y$ we define $\|y(t)\|_Y = \sup_{t\in[t_0,\infty)} \|y(t)\|_{DL_p(R^3)}$. Let $\{y_n|n = 1.2.\dots\} \subset Y$be a Cauchy sequence in $Y$. For every $t \in [t_0,\infty)$, $\{y_n(t)|n = 1,2,\dots\}$ is a Cauchy sequence in $DL_p(R^3)$, so there exists $y_0(t) \in DL_p(R^3)$ such that $\lim_{n\to\infty} y_n(t) = y_0(t)$ uniformly on $[0,\infty)$, i.e., $\lim_{n\to\infty}\|y_n(t) - y_0(t)\|_{DL_p(R^3)} = 0$ uniformly on $[t_0,\infty)$. We now prove that $y_0(t) \in Y$.

Because $\lim_{n\to\infty}\|y_n(t)\|_{DL_p(R^3)} = \|y_0(t)\|_{DL_p(R^3)}$ is uniform on $[t_0,\infty)$, using Theorem7.11 in [20] we have

$$\lim_{t\to t_1}\|y_0(t)\|_{DL_p(R^3)} = \lim_{t\to t_1}\left\|\lim_{n\to\infty} y_n(t)\right\|_{DL_p(R^3)} = \lim_{t\to t_1}\lim_{n\to\infty}\|y_n(t)\|_{DL_p(R^3)}$$

$$= \lim_{n\to\infty}\lim_{t\to t_1}\|y_n(t)\|_{DL_p(R^3)} = \lim_{n\to\infty}\|y_n(t_1)\|_{DL_p(R^3)} = \|y_0(t_1)\|_{DL_p(R^3)}.$$

For any continuous functional $\phi \in (DL_p(R^3))^*$, we have

$$\lim_{t\to t_1}\phi\big(y_0(t)\big) = \lim_{t\to t_1}\phi\left(\lim_{n\to\infty} y_n(t)\right) = \lim_{t\to t_1}\lim_{n\to\infty}\phi\big(y_n(t)\big) = \lim_{n\to\infty}\lim_{t\to t_1}\phi\big(Y_n(t)\big)$$

$$= \lim_{n\to\infty}\phi(\lim_{t\to t_1}(y_n(t)) = \lim_{n\to\infty}\phi\big(y_n(t_1)\big)$$

$$= \phi\left(\lim_{n\to\infty} y_n(t_1)\right) = \phi(y_0(t_1)).$$

So $y_0(t)$ weakly converges to $y_0(t_0)$. It follows from Proposition 1.17 in [4] that

$$\lim_{t\to t_1}\|y_0(t) - y_0(t_1)\|_{DL_p(R^3)} = 0.$$

Hence, $y_0(t)$ is continuous.

(2) It follows from Theorems 7.16 and 7.11 in [21] that

$$\lim_{t\to\infty}\int_{t_0}^{t}\|y_0(s)\|_{DL_p(R^3)}ds$$

$$= \lim_{t\to\infty}\int_{t_0}^{t}\left\|\lim_{n\to\infty} y_n(s)\right\|_{DL_p(R^3)}ds$$

$$= \lim_{t\to\infty}\int_{t_0}^{t}\lim_{n\to\infty}\|y_n(s)\|_{DL_p(R^3)}ds$$

$$= \lim_{t\to\infty}\lim_{n\to\infty}\int_{t_0}^{t}\|y_n(s)\|_{DL_p(R^3)}ds$$

$$= \lim_{n\to\infty}\lim_{t\to\infty}\int_{t_0}^{t}\|y_n(s)\|_{DL_p(R^3)}ds,$$

so $\int_{t_0}^{\infty}\|y_0(t)\|_{DL_p(R^3)}\ dt$ converges. Hence, $y_0(t) \in Y$, and therefore $Y$ is a Banach space.

Now, we prove that $f(t, A^{-\alpha}y(t)) \in Y$ for $y(t) \in Y$, where

$$f(t, A^{-\alpha}y(t)) = -(A^{-\alpha}y(t) \bullet \nabla)A^{-\alpha}y(t).$$

First, we prove that $A^{-\alpha}y(t) \in Y$. Since $A^{-\alpha}$ is bounded and continuous, $A^{-\alpha}y(t)$ is continuous. From (12), $f(t, A^{-\alpha}y(t))$ is also continuous for $t \in [t_0,\infty)$, and

$$\lim_{t\to\infty}\int_{t_0}^{t}\|A^{-\alpha}y(s)\|_{DL_p(R^3)}ds$$

$$\leq \lim_{t\to\infty}\int_{t_0}^{t}\|A^{-\alpha}\|_{B(DL_p(R^3))}\|y(s)\|_{DL_p(R^3)}ds$$

$$= \|A^{-\alpha}\|_{B(DL_p(R^3))}\lim_{t\to\infty}\int_{t_0}^{t}\|y(s)\|_{DL_p(R^3)}ds.$$

Hence $\int_{t_0}^{\infty}\|A^{-\alpha}y(s)\|_{DL_p(R^3)}ds$ converges, so $A^{-\alpha}y(t)\in Y$. $\int_{t_0}^{\infty}\|y(s)\|_{DL_p(R^3)}ds$ is convergent， $\|y(t)\|_{DL_p(R^3)}$ is bounded on $[t_0,\infty)$, i.e., there exists a constant $K$ such that $\|y(t)\|_{DL_p(R^3)}\leq K$ for $t\in[t_0,\infty)$. Therefore, we have

$$\lim_{t\to\infty}\int_{t_0}^{t}\left\|f\left(s,A^{-\alpha}y(s)\right)\right\|_{DL_p(R^3)}ds$$

$$= \lim_{t\to\infty}\int_{t_0}^{t}\|(A^{-\alpha}y(s)\bullet\nabla)A^{-\alpha}y(s)\|_{DL_p(R^3)}ds$$

$$\leq \lim_{t\to\infty}\int_{t_0}^{t}M\left\|A^{-\frac{1}{2}}y(s)\right\|^2_{DL_p(R^3)_{\frac{1}{2}}}ds$$

$$= M\lim_{t\to\infty}\int_{t_0}^{t}\|y(s)\|^2_{DL_p(R^3)}ds$$

$$\leq MK\lim_{t\to\infty}\left[\int_{t_0}^{t}\|y(s)\|^2_{DL_p(R^3)}ds\right]$$

$$= MK\int_{t_0}^{\infty}\|y(s)\|^2_{DL_p(R^3)}ds.$$

Hence, $\int_{t_0}^{\infty}\|f(s,A^{-\alpha}y(s)\|_{DL_p(R^3)}ds$ converges. Therefore, $f(t,A^{-\alpha}y(t))\in Y$.

(Step 3) On $Y$ we define a mapping $F$ by

$$Fy(t)=\begin{cases}\frac{1}{t+1}T(t)A^{\alpha}u_0+\frac{1}{t+1}\int_{t_0}^{t}A^{\alpha}T(t-s)f(s,A^{-\alpha}y(s))ds, t\in[t_0,t_1]\\ \frac{1}{t+1}T(t)A^{\alpha}u_0+\frac{1}{t+1}\int_{t_0}^{t_1}A^{\alpha}T(t-s)f(s,A^{-\alpha}y(s))ds, t\geq t_1.\end{cases} \quad (18)$$

We prove that $F:Y\to Y$. $Fy(t)$ is continuous for $t\in[t_0,\infty)$. Since $f(t,A^{-\alpha}y(t))\in Y$, from Theorem 2.6.13 (c) in [18], we see that $A^{\alpha}T(t-s)$ is bounded and

$$\|A^{\alpha}T(t-s)\|_{B(DL_p(R^3)}\leq M_{\alpha}(t-s)^{-\alpha}e^{-\delta(t-s)},$$

where $M_{\alpha}$ is a constant and $\delta>0$. Therefore

$$\frac{1}{s+1}\|A^{\alpha}T(s-s')\|_{B\left(DL_p(R^3)\right)}\leq\frac{(s-s')^{-\alpha}}{s+1}M_{\alpha}e^{-\delta(s-s')}.$$

Since $\alpha=\frac{1}{2}$, when $t\to\infty (s\in[t_0,t_1])$ $\frac{(t-s)^{-\alpha}}{t+1}\to 0$ and $e^{-\delta(t-s)}<1$. Hence, $\frac{1}{t+1}$ $\|A^\alpha T(t-s)\|_{B(DL_p(R^3))}$ is bounded, and $\frac{1}{t+1}\|A^\alpha T(t-s)\|_{DL_p(R^3)}\le K$ for some constant $K$. Therefore,

$$\lim_{t\to\infty}\int_{t_0}^{t}\|Fy(s)\|_{DL_p(R^3)}ds$$

$$\le \lim_{t\to\infty}\frac{1}{s+1}\int_{t_0}^{t}\|T(s)A^\alpha u_0\|_{DL_p(R^3)}\,ds+\lim_{t\to\infty}\int_{t_0}^{t}(\int_{t_0}^{t_1}\frac{1}{s+1}\|A^\alpha T(s-s')f(s',A^{-\alpha}y(s')\|_{DL_p(R^3)}ds')d$$

$$\le \lim_{t\to\infty}\int_{t_0}^{t}\frac{1}{s+1}\|T(s)\|_{B(DL_p(R^3))}\|A^\alpha u_0\|_{DL_p(R^3)}ds+\lim_{t\to\infty}\int_{t_0}^{t}(\int_{t_0}^{t_1}\frac{1}{s+1}\|A^\alpha T(s-s')f(s',A^{-\alpha}y(s')\|_{DL_p(R^3)}ds')ds$$

$$=\lim_{t\to\infty}\int_{t_0}^{t}(\int_{t_0}^{t_1}\frac{1}{s+1}\|A^\alpha T(s-s')\|_{B(DL_p(R^3))}\|f(s',A^{-\alpha}y(s')\|_{DL_p(R^3)}ds')ds$$

$$\le K\lim_{t\to\infty}\int_{t_0}^{t}(\int_{t_0}^{t_1}\|f(s',A^{-\alpha}y(s')\|_{DL_p(R^3)}ds')ds.$$

Hence, $\int_{t_0}^{\infty}\|Fy(s)\|_{DL_p(R^3)}\,ds$ converges, and therefore $Fy(t)\in F$.

For every $y\in Y$, $Fy(0)=A^\alpha u_0$. Let $S$ be the nonempty closed and bounded subset of $Y$ defined by

$$S=\{y: y\in Y, y(0)=A^\alpha u_0, \|y(t)-A^\alpha u_0\|_Y\le\delta\}. \quad (19)$$

For $y\in S$, when $t\in[t_0,t_1]$, using (12), (14), (17) and (19) we have

$$\|Fy(t)-A^\alpha u_0\|_Y$$

$$\le\left\|\frac{1}{t+1}T(t)A^\alpha u_0-A^\alpha u_0\right\|_Y+\left\|\int_{t_0}^{t}A^\alpha T(t-s)[f(s,A^{-\alpha}y(s)-f(s,u_0)]ds\right\|_Y+\left\|\int_{t_0}^{t}A^\alpha T(t-s)f(s,u_0)ds\right\|_Y$$

$$\le\frac{\delta}{2}+C_\alpha(L\delta+B)\int_{t_0}^{t}(t-s)^{-\alpha}\,ds=\frac{\delta}{2}+C_\alpha(1-\alpha)^{-1}(L\delta+B)t_1^{1-\alpha}\le\delta.$$

For $t\in[t_1,\infty)$, we have

$$\|Fy(t)-A^\alpha u_0\|_Y$$

$$\le\left\|\frac{1}{t+1}T(t)A^\alpha u_0-A^\alpha u_0\right\|_Y+\left\|\int_{t_0}^{t_1}A^\alpha T(t-s)[f(s,A^{-\alpha}y(s)-f(s,u_0)]ds\right\|_Y+\left\|\int_{t_0}^{t_1}A^\alpha T(t-s)f(s,u_0)ds\right\|_Y$$

$$\le\frac{\delta}{2}+C_\alpha(L\delta+B)\int_{t_0}^{t_1}(t-s)^{-\alpha}\,ds=\frac{\delta}{2}+C_\alpha(1-\alpha)^{-1}(L\delta+B)t_1^{1-\alpha}\le\delta.$$

Therefore $F:S\to S$. Furthermore, if $y_1,y_2\in S$ and $t\in[t_0,t_1]$ then

$$\|Fy_1(t)-Fy_2(t)\|_Y$$

$$\le\int_{t_0}^{t}\|A^\alpha T(t-s)\|_Y\cdot\|f(s,A^{-\alpha}y_1(s))-f(s,A^{-\alpha}y_2(s))\|_Y ds$$

$$\le LC_\alpha(1-\alpha)^{-1}{t_1}^{1-\alpha}\|y_1-y_2\|_Y\le\frac{1}{2}\|y_1-y_2\|_Y. \quad (20)$$

$$\|Fy_1(t)-Fy_2(t)\|_Y$$

$$\leq \int_{t_0}^{t_1}\|A^{\alpha}T(t-s)\|_{\ Y}\cdot\left\|f\big(s,A^{-\alpha}y_1(s)\big)-f\big(s,A^{-\alpha}y_2(s)\big)\right\|_Y \mathrm{ds}$$

$$\leq \frac{1}{2}\|y_1-y_2\|_Y.$$

(see $p.197^1-197^{13}$in [18]) From the contraction mapping theorem, we see that the mapping $F$ has a unique fixed point $y_0(t)\in S$ and $y_0(t)$ satisfies the integral equation

$$y_0(t)=\begin{cases} y_1(t)=\frac{1}{t+1}T(t)A^{\alpha}u_0+\frac{1}{t+1}\int_{t_0}^{t}A^{\alpha}T(t-s)f\big(s,A^{-\alpha}y_0(s)\big)ds, t\in[t_0,t_1] \\ y_2(t)=\frac{1}{t+1}T(t)A^{\alpha}u_0+\frac{1}{t+1}\int_{t_0}^{t_1}A^{\alpha}T(t-s)f\big(s,A^{-\alpha}y_0(s)\big)ds, t\geq t_1. \end{cases} \tag{21}$$

Since $y_0(t)\in Y$, then $f(t,A^{-\alpha}y_0(t))\in Y$. From (12), and because $y$ is continuous, it follows that $t\to f(t,A^{-\alpha}y(t))$ is continuous on $[t_0,t_1]$ and bounded on $[t_0,t_1]$. Therefore we can let

$$\|f(t,A^{-\alpha}y_0(t)\|_Y\leq N \quad \text{for } t\in[t_0,t_1], \tag{22}$$

where $N$ is a constant.

(Step 4) We can prove that the solution $y_0(t)$ is locally Hölder continuous about $t$ on $(t_0,\infty)$. By Theorem 2.6.13 (d) and (c) in [18] for every $\beta$ satisfying $0<\beta<1-\alpha$ and every $0<h<1$ we have

$$\|(T(h)-I)A^{\alpha}T(t-s)\|_Y\leq C_\beta h^\beta\left\|A^{\alpha+\beta}T(t-s)\right\|_Y\leq Ch^\beta(t-s)^{-(\alpha+\beta)}. \tag{23}$$

If $t_0< t < t+h\leq t_1$, then

$$\|Y_0(t+h)-Y_0(t)\|_Y$$

$$\leq \|\ (\frac{1}{t+h+1}T(t+h)-\frac{1}{t+h+1}T(t)+\frac{1}{t+h+1}T(t)-\frac{1}{t+1}\mathrm{T(t)})A^{\alpha}u_0\ \|\ ds$$

$$+\int_{t_0}^{t}\|(T(h)-I)A^{\alpha}T(t-s)f(s,A^{-\alpha}y_0(s))\|\ ds$$

$$+\int_{t}^{t+h}\|A^{\alpha}T(t+h-s)f(s,A^{-\alpha}y_0(s)\|$$

$$\leqslant\frac{h}{(t+h+1)(t+1)}\|T(t)A^{\alpha}u_0\|_Y+\frac{1}{t+h+1}I_1+I_2+I_3$$

$$=I_0+I_1+I_2+I_3. \tag{24}$$

Using (22) and (23) we estimate each term of (24).

$$I_0\leq\|T(t)\|_{B(DL_p(R^3))}\|A^{\alpha}u_0\|_{DL_p(R^3)}h\leq M_0h^\beta, \tag{25}$$

$$I_1\leq Ct^{-(\alpha+\beta)}\|u_0\|_{DL_p(R^3)}h^\beta\leq M_1h^\beta, \tag{26}$$

$$I_2\leq CN\int_{t_0}^{t'}(t-s)^{-(\alpha+\beta)}\,ds\leq M_2h^\beta, \tag{27}$$

$$I_3 \leq NC_\alpha \int_{t'}^{h'} (t \ \ + h - s)^{-\alpha}\, ds = \frac{NC_\alpha}{1-\alpha} h^{1-\alpha} \leq \ M_3 h^\beta. \tag{28}$$

Combining (24)-(28) it follows that for every $t' > t_0$, $t' < t, s \leq t_1$ there is a constant $C$ such that

$$\|y_0(t) - y_0(s)\|_Y \leq C|t-s|^\beta \qquad \text{for } t' < t, s \leq t_1. \tag{29}$$

For $t \in (t_0, t_1], s \geq t_1$ and $t, s \geq t_1$ we also have the formula (27). Therefore, $y_0$ is local Hölder continuous on $(t_0, \infty)$. (see p.198 in [18])

From (12) and (29), we have

$$\begin{aligned} &\|f(s, A^{-\alpha} y_0(s) - f(t, A^{-\alpha} y_0(t)\|_Y \\ &\leq L(|t-s|^\vartheta + \|y_0(t) - y_0(s)\|_Y) \\ &\leq C_1(|t-s|^\vartheta + |t-s|^\beta). \end{aligned}$$

Hence $t \to f(t, A^{-\alpha} y_0(t))$ is local Hölder continuous. Let $y_0(t)$ be the solution of (21). The Corollary 4.3.3 in [18] is also valid for $f \in L^1((0,\infty): X)$. Thus, the inhomogeneous initial problem

$$\begin{cases} \dfrac{du(t)}{dt} + Au(t) = f(t, A^{-\alpha} y_0(t)), \\ \qquad\qquad u(t_0) = u_0 \end{cases}$$

has a unique solution

$$u \in C^1((t_0,\infty): DL_p(R^3)), \text{ and } u(t) = T(t-t_0)u_0 + \int_{t_0}^{t} T(t-s) f\big(s, A^{-\alpha} y_0(s)\big) ds. \tag{30}$$

Therefore, we have

$$A^\alpha u(t) = T(t-t_0) A^\alpha u_0 + \int_{t_0}^{t} A^\alpha T(t-s) f(s, A^{-\alpha}\, y_0(s)) ds$$

$$= (t+1) y_1(t), t \in [t_0, t_1],$$

$$A^\alpha u(t) = T(t-t_0) A^\alpha u_0 + \int_{t_0}^{t_1} A^\alpha T(t-s) f(s, A^{-\alpha}\, y_0(s)) d$$

$$= (t+1) y_2(t), t \in [t_1, \infty), \ (31)$$

From this, we have $A^\alpha u(t) = y_0(t)$, so

$$u(t) = \begin{cases} (t+1) A^{-\alpha} y_1(t), t \in [t_0, t_1] \\ (t+1) A^{-\alpha} y_2(t), t \in [t_1, \infty) \end{cases} \in C([t_0,\infty): DL_p(R^3)) \cap C^1((t_0,\infty): DL_p(R^3))$$

is a global solution of (13). ☐

Below, we need the Banach lattice (see [3] [16]). A real vector space $G$ that is ordered by some order relation $\leq$ is called a vector lattice (or Riesz space) if any two elements $f, g \in G$ have a least upper bound, denoted by $f \vee g$, and a greatest lower bound, denoted by $f \wedge g$, and the following properties are satisfied:

(i) If $f \leq g$, then $f + h \leq g + h$ for all $f, g, h \in G$,

(ii) If $0 \leq f$, then $0 \leq tf$ for all $f \in G$ and $0 \leq t \in R$.

**Definition 3.** A Banach lattice is a real Banach space $G$ endowed with an ordering $\le$ such that $(G, \le)$ is a vector lattice and the norm is a lattice norm, i.e., $|f| \le |g|$ implies $\|f\| \le \|g\|$ for $f, g \in G$, where $|f| = f \vee (-f)$ is the absolute value of $f$ and $\|\bullet\|$ is the norm in $G$.

In a Banach lattice $G$ we define for $f \in G$

$$f^+ := f \vee 0, \qquad f^- := (-f) \vee 0.$$

The absolute value of $f$ is $|f| = f^+ + f^-$ and $f = f^+ - f^-$.

$$0 \le f \le g \Rightarrow |f| \le |g| \Rightarrow \|f\| \le \|g\|.$$

We will need the above formula below. In $L^p(R^3)$ we define the order for $f, g \in L^p(R^3)$,

$$f \le g \Leftrightarrow f(x) \le g(x) \text{ for } a.e.\ x \in R^3,$$

and $(f \vee g)(x) := max\{f(x), g(x)\}$, $(f \wedge g)(x) := min\{f(x), g(x)\}$ for $a.e.$ $x \in R^3$. Then $L^p(R^3)$ and $L_p(R^3)$ are all Banach lattices (see p.148 in [3]).

**Lemma 5.** Suppose that $u, v \in DL_p(R^3)$ are smooth, then $(u \bullet \nabla)v \in DL_p(R^3)$.

Proof. Suppose that $u, v \in DL_p(R^3)$ are smooth. From Theorem3(a) in [7, p.129], we have

$$\int_{R^3} u \bullet \nabla h dx = 0, \int_{R^3} v \bullet \nabla h dx = 0$$

for all $h \in W^{1,p}(R^3)$. That is, $\int_{R^3} \sum_{i=1}^{3} u_i \frac{\partial h}{\partial x_i} dx = 0$. Since $L_p(R^3)$ is a Banach lattice and $\nabla h \in DL_p(R^3)^{\perp}$ for $h \in W^{1,p}(R^3)$, then $u = u^+ - u^-, u^+, u^- \in DL_p(R^3), \nabla h = (\nabla h)^+ - (\nabla h)^-$. From Lemma 7.6 in [8], we have

$$(\nabla h)^+ = \begin{pmatrix} \frac{\partial h}{\partial x_1} \\ \frac{\partial h}{\partial x_2} \\ \frac{\partial h}{\partial x_3} \end{pmatrix}^+ = \begin{pmatrix} (\frac{\partial h}{\partial x_1})^+ \\ (\frac{\partial h}{\partial x_2})^+ \\ (\frac{\partial h}{\partial x_3})^+ \end{pmatrix} = \begin{pmatrix} \frac{\partial h^+}{\partial x_1} \\ \frac{\partial h^+}{\partial x_2} \\ \frac{\partial h^+}{\partial x_3} \end{pmatrix} = \nabla(h^+).$$

Similar $(\nabla h)^- = \nabla(h^-)$. $W^{1,p}(R^3) \subset W^{0,p}(R^3) = L_p(R^3)$ are all Banach spaces. $h \in W^{1,p}(R^3)$ implies

$h^+, h^- \in W^{1,p}(R^3), \int_{R^3} u^+ \bullet \nabla(h^+) dx = 0.$

Let $\Omega_n (n = 1,2,3, \ldots)$ be a bounded closed ball with a radius of $n$, then $\Omega_n \subset \Omega_{n+1} (n = 1,2,3, \ldots.)$, $\bigcup_1^{\infty} \Omega_n = R^3$, and $\Omega_n \to R^3$ when $n \to \infty$. $v$ *is smooth*, $\frac{\partial v_i}{\partial x_j}$ are all continuous on $\Omega_n (i, j = 1,2,3, n = 1,2,3, \ldots)$. Hence $\frac{\partial v_i}{\partial x_j}$ are bounded on $\Omega_n (i = 1,2,3, n = 1,2,3, \ldots)$. Thus $\left|\frac{\partial v_i}{\partial x_{i,j}}\right| \le L_1$ $(i, j = 1,2,3)$ for some constant $L_1 > 0$. Let

$$L_2 = Max\left\{ \sup_{x \in \Omega_n} u_1^+, \sup_{x \in \Omega_n} u_2^+, \sup_{x \in \Omega_n} u_3^+ \right\}.$$

We then have

$$0=\int_{\Omega_n} u^{+}\cdot\nabla(h^{+})\,dx$$

$$=-L_1\int_{\Omega_n}\sum_{i=1}^{3} u_i{}^{+}\,(\frac{\partial h}{\partial x_i})^{+}dx$$

$$\leq\int_{\Omega_n}\sum_{i=1}^{3} u_i{}^{+}\left|\frac{\partial v_i}{\partial x_i}\right|\ \ (\frac{\partial h}{\partial x_i})^{+}dx$$

$$\leq\int_{\Omega_n}\sum_{i=1}^{3}(\sum_{j=1}^{3} u_j{}^{+}\left|\frac{\partial v_i}{\partial x_j}\right|)\,(\frac{\partial h}{\partial x_i})^{+}dx$$

$$\leq 3L_1L_2\int_{\Omega_n}\sum_{i=1}^{3} u_i{}^{+}\left|\frac{\partial v_i}{\partial x_i}\right|(\frac{\partial h}{\partial x_i})^{+}dx$$

$$\leq 3L_1^2L_2\int_{\Omega_n}\sum_{i=1}^{3} u_i{}^{+}(\frac{\partial h}{\partial x_i})^{+}dx=0$$

It follows that

$$\int_{\Omega_n}\left|\sum_{i=1}^{3}(\sum_{j=1}^{3} u_j{}^{+}\frac{\partial v_i}{\partial x_j})(\frac{\partial h}{\partial x_i})^{+}\right|dx$$

$$\leq\int_{\Omega_n}\sum_{i=1}^{3}(\sum_{j=1}^{2} u_j{}^{+}\left|\frac{\partial v_i}{\partial x_j}\right|)(\frac{\partial h}{\partial x_i})^{+}dx=0.$$

Hence

$\int_{\Omega_n}(u^{+}\bullet\nabla)v\bullet(\nabla h)^{+}dx=0$(n=1, 2, 3, ...) and

$$\int_{R^3}(u^{+}\bullet\nabla)v\bullet(\nabla h)^{+}dx$$

$$=\lim_{n\to\infty}\int_{\Omega_n}(u^{+}\bullet\nabla)v\bullet(\nabla h)^{+}dx=0.$$

Similarly, $\int_{R^3}(u^{+}\bullet\nabla)v\bullet(\nabla h)^{-}dx=0$. Then we have $\int_{R^3}(u^{+}\bullet\nabla)v\bullet\nabla h dx=\int_{R^3}(u\bullet\nabla)v\bullet[(\nabla h)^{+}-(\nabla h)^{-}]\,dx=0$. Similarly, $\int_{R^3}(u^{-}\bullet\nabla)v\bullet\nabla h dx=0$. Therefore,

$$\int_{R^3}(u\bullet\nabla)v\bullet\nabla h dx=\int_{R^3}((u^{+}-u^{-})\bullet\nabla)v\bullet\nabla h dx$$

$$=\int_{R^3}(u^{+}\bullet\nabla)v\bullet\nabla h dx-\int_{R^3}(u^{-}\bullet\nabla)v\bullet\nabla h dx=0,$$

and so $(u\bullet\nabla)v\in DL_p(R^3)$. ☐

## III. MAIN RESULT

We now consider the condition for the N-S equation to have a unique smooth global solution is the initial value $(t_0, u_0) \in U$. Let

$$H\left([t_0, \infty); C^\infty(R^3) \cap DL_p(R^3)_{1/2}\right)$$

denote the space of all Hölder continuous functions $u(t)$ on $[t_0, \infty)$ with different exponents in $(0,1]$ and with smooth function values in the Banach space $DL_p(R^3)_{1/2}$. From formula (2), we can see that these function values are all divergence-free. Then from Lemma 5, for any $u \in H\left([t_0, \infty); C^\infty(R^3) \cap DL_p(R^3)_{1/2}\right)$ and any $t_1, t_2 \in [0, \infty)$, $(u(t_1) \bullet \nabla)u(t_2) \in DL_p(R^3)$; and for any $u_1, u_2 \in H\left([t_0, \infty); C^\infty(R^3) \cap DL_p(R^3)_{1/2}\right)$ and any $t \in [0, \infty)$, $(u_1(t) \bullet \nabla)u_2(t) \in DL_p(R^3)$. Thus, the bilinear form $(v(t) \bullet \nabla)u(t)$ on $DL_p(R^3)_{1/2}$ takes values in $DL_p(R^3)$.

Suppose that $u(x) = (u_1(x_1), u_2(x_2), u_3(x_3)) \in DL_p(R^3)_{1/2}$ is smooth. Let $u(t,x) \equiv u(x)$ for $t \in [t_0, \infty)$. Then $u(t,x) \in H\left([t_0, \infty); C^\infty(R^3) \cap DL_p(R^3)_{1/2}\right)$. Hence $H\left([t_0, \infty); C^\infty(R^3) \cap DL_p(R^3)_{1/2}\right)$ is not empty and has nonzero functions. Consider the graph

$$U := \left\{(t, u(t)) : t \in (t_0, \infty), \forall u \in H\left([t_0, \infty); C^\infty(R^3) \cap DL_p(R^3)_{1/2}\right)\right\}.$$

$U$ is a subset of $(t_0, \infty) \times DL_p(R^3)_{1/2}$. For any smooth function $u(x) = \left(u_1(x_1), u_2(x_2), u_3(x_3)\right) \in DL_p(R^3)_{\frac{1}{2}}$, we let $u(t,x) \equiv u(x)$ for all $t \in [t_0, \infty)$. Then $u(t) \in H\left([t_0, \infty); C^\infty(R^3) \cap DL_p(R^3)_{1/2}\right)$ and $(t_0, u) \in U$. Hence, $U$ is not empty and we have nonzero functions $F\left(t, u(t)\right) = -(u(t) \bullet \nabla)u(t)$ from $U$ to $DL_p(R^3)_{1/2}$.

**Theorem 1.** Suppose that the given externally applied force $f(t,x) = 0$. In the framework function space $L_p(R^3)$, the Navier-Stokes equation (1) has a unique smooth global solution $p(t,x), u(t,x)$ on $[t_0, \infty) \times R^3$ with bounded energy if and only if the initial value $u_0$ is smooth and divergence free vector field.

Proof. Sufficiency. Suppose that the initial value $u_0$ is smooth, divergence free vector field.

(Step 1) First, $F(t, u(t)) = -(u(t) \bullet \nabla)u(t)$ is a function: $U \to DL_p(R^3)$ because $(u(t) \bullet \nabla)u(t) \in DL_p(R^3)$, according to the definition of $U$.

We find that by incorporating the divergence-free condition, we can remove the pressure term from our equation. (see p. $271^3$ in [5] and pp. $234_6$ and $239_9$ in [19]). From $DL_p(R^3)^\perp = \{\nabla h; h \in W^{1,2}(R^3)\}$, we see that $\nabla p \in DL_p(R^3)^\perp$, and so $P\nabla p = 0$. For $u \in DL_p(R^3)$, we have $\Delta u \in DL_p(R^3)$ because of Lemma 4. Hence, by applying $P$ to equation (1) we have $P\Delta u = \Delta u$. It follows from definition of $U$ that $(u \bullet \nabla)u \in DL_p(R^3)$, so $P(u \bullet \nabla)u = (u \bullet \nabla)u$. Therefore we can rewrite (1) as an abstract initial value problem on $DL_p(R^3)$,

$$\begin{cases}\frac{du}{dt} = \Delta u + F(t,u(t)), t \in (t_0,\infty), \\ \qquad u|_{t=t_0} = u_0, x \in R^3\end{cases} \tag{32}$$

where $F(t,u(t)) = -(u(t) \bullet \nabla)u(t)$. From Lemma 2, $\Delta_{|DL_p(R^3)}$ is the generator of an analytic semigroup $T(t)$ of contraction on $DL_p(R^3)$. Therefore, $\|T(t)\| \le 1$.

(Step 2) If $u \in H\left([t_0,\infty); C^\infty(R^3) \cap DL_p(R^3)_{1/2}\right)$, then there is a constant $C$ and 0<$\beta \le 1$ such that

$$\|u(t_1,x) - u(t_2,x)\|_{DL_p(R^3)_{1/2}}$$

$$\le C|t_1 - t_2|^\beta \text{ for } t_1,t_2 \in [0,\infty). \tag{33}$$

For any$(t_1,u_1(t_1)),(t_2,u_2(t_2)) \in U$we have $(u_1(t_1) \bullet \nabla)u_1(t_1), (u_2(t_1) \bullet \nabla)u_2(t_1) \in DL_p(R^3)$, and

$$\|(u_1(t_1) \bullet \nabla)u_1(t_1) - (u_2(t_1) \bullet \nabla)u_2(t_1)\|_{DL_p(R^3)}$$
$$= \|(u_1(t_1) \bullet \nabla)u_1(t_1) - (u_1(t_1) \bullet \nabla)u_2(t_1)\|_{DL_p(R^3)}$$
$$+\|(u_1(t_1) \bullet \nabla)u_2(t_1) - (u_2(t_1) \bullet \nabla)u_2(t_1)\|_{DL_p(R^3)}$$
$$\le \left\|(u_1(t_1) \bullet \nabla)\left(u_1(t_1) - u_2(t_1)\right)\right\|_{DL_p(R^3)}$$
$$+\|[(u_1(t_1) - u_2(t_1)) \bullet \nabla)]u_2(t_1)\|_{DL_p(R^3)}$$
$$\le M(\|u_1(t_1)\|_{DL_p(R^3)_{1/2}} \left\|\left(u_1(t_1) - u_2(t_1)\right)\right\|_{DL_p(R^3)_{1/2}}$$
$$+\|u_1(t_1) - u_2(t_1)\|_{DL_p(R^3)_{1/2}} \|u_2(t_1)\|_{DL_p(R^3)_{1/2}})$$
$$= M\left(\|u_1(t_1)\|_{DL_p(R^3)_{1/2}} + \|u_2(t_1)\|_{DL_p(R^3)_{1/2}}\right)$$
$$\|u_1(t_1) - u_2(t_1)\|_{DL_p(R^3)_{1/2}}$$
$$\le M\left(\|u_1(t_1)\|_{DL_p(R^3)_{1/2}} + \|u_2(t_2)\|_{DL_p(R^3)_{1/2}} + \right.$$
$$\|-u_2(t_2) + u_2(t_1)\|_{DL_p(R^3)_{1/2}})$$
$$\|u_1(t_1) - u_2(t_1)\|_{DL_p(R^3)_{1/2}}$$
$$\le M\left(\|u_1(t_1)\|_{DL_p(R^3)_{1/2}} + \|u_2(t_2)\|_{DL_p(R^3)_{1/2}} + C|t_1 - t_2|^{\beta_1}\right)\|u_1(t_1) - u_2(t_1)\|_{DL_p(R^3)_{1/2}}. \tag{34}$$

We used Lemma 3 in the above third step and formula (33) in the sixth step.

For any $u \in H\left([0,\infty); C^\infty(R^3) \cap DL_p(R^3)_{1/2}\right)$ we have

$$\|(u(t_1) \bullet \nabla)u(t_1) - (u(t_2) \bullet \nabla)u(t_2)\|_{DL_p(R^3)}$$
$$\le \|(u(t_1) \bullet \nabla)u(t_1) - (u(t_1) \bullet \nabla)u(t_2)\|_{DL_p(R^3)}$$
$$+\|(u(t_1) \bullet \nabla)u(t_2) - (u(t_2) \bullet \nabla)u(t_2)\|_{DL_p(R^3)}$$

$$= \left\|\left(u(t_1) \bullet \nabla\right)\left(u(t_1) - u(t_2)\right)\right\|_{DL_p(R^3)}$$

$$+\left\|\left[\left(u(t_1) - u(t_2)\right) \bullet \nabla\right]u(t_2)\right\|_{DL_p(R^3)}$$

$$\leq M\left(\|u(t_1)\|_{DL_p(R^3)_{1/2}} \|u(t_1) - u(t_2)\|_{DL_p(R^3)_{1/2}}\right.$$

$$\left.+\|u(t_1) - u(t_2)\|_{DL_p(R^3)_{1/2}} \|u(t_2)\|_{DL_p(R^3)_{1/2}}\right)$$

$$= M(\|u(t_1)\|_{DL_p(R^3)_{1/2}} + \|u(t_2)\|_{DL_p(R^3)_{1/2}})$$

$$\|u(t_1) - u(t_2)\|_{DL_p(R^3)_{1/2}}$$

$$\leq MC(\|u(t_1)\|_{DL_p(R^3)_{1/2}} + \|u(t_2)\|_{DL_p(R^3)_{1/2}})|t_1 - t_2|^{\beta}. \quad (35)$$

We used Lemma 3 in above third step and formula (33) in the fifth step.

(Step 3) If the initial value $u_0$ is smooth and divergence free vector field, then $u_0 \in C_{0,\sigma}^{\infty}(R^3)$. *From the formula* (10)*we see that* $u_0 \in DL_p(R^3)_{1/2}$.

From the above we see that if let $u(t,x) \equiv u_0$ for all $t \in [t_0,\infty)$, then $u(t,x) \in H\left([t_0,\infty); C^{\infty}(R^3) \cap DL_p(R^3)_{1/2}\right)$ and $(t_0,u_0) \in U$. Set

$$V = B_{\varepsilon}(t_0,u_0) = \{(t,u(t)) \in U: |t - t_0\ | < \varepsilon \leq 1, \|u - u_0\|_{DL_p(R^3)_{1/2}} < \varepsilon\}.$$

Then for $(t,u(t)) \in V$,

$$\|u\|_{DL_p(R^3)_{1/2}} = \|u - u_0 + u_0\|_{DL_p(R^3)_{1/2}}$$

$$\leq \|u - u_0\|_{DL_p(R^3)_{1/2}} + \|u_0\|_{DL_p(R^3)_{1/2}} \leq \varepsilon + \|u_0\|_{DL_p(R^3)_{1/2}}.$$

Let $L = \varepsilon + \|u_0\|_{DL_p(R^3)_{1/2}}$, then formulas（34）and (35) become

$$\|(u_1(t_1) \bullet \nabla)u_1(t_1) - (u_2(t_1) \bullet \nabla)u_2(t_1)\|_{DL_p(R^3)}$$

$$\leq M(2L + C_1|t_1 - t_2|^{\beta_1})\|u_1(t_1) - u_2(t_1)\|_{DL_p(R^3)_{1/2}}$$

$$\left\|\ (u(t_1) \bullet \nabla)u(t_1) - (u(t_2) \bullet \nabla)u(t_2)\right\|_{DL_p(R^3)}$$

$$\leq M2LC|t_1 - t_2|^{\beta}. \quad (36)$$

Let

$$L_1 = M(2L + 2^{\beta_1}),$$

$$L_2 = M(2L + 2^{\beta_1})C_2 + 2MLC,$$

$$\beta_3 = min(\beta, \beta_2), L_3 = Max(L_1, L_2).$$

From (36), and because $|t_1 - t_2| \leq |t_1 - t_0| + |t_0 - t_2| \leq 2\varepsilon \leq 2$ for all $(t_i, u_i) \in V (i = 1,2)$, we have

$$\begin{aligned}
&\|F(t_1, u_1(t_1)) - F(t_2, u_2(t_2))\|_{DL_p(R^3)} \\
&\leq \|F(t_1, u_1(t_1)) - F(t_1, u_2(t_1))\|_{DL_p(R^3)} \\
&\quad + \left\|F\left(t_1, u_2(t_1)\right) - F\left(t_2, u_2(t_2)\right)\right\|_{DL_p(R^3)} \\
&= \|(u_1(t_1) \bullet \nabla)u_1(t_1) - (u_2(t_1) \bullet \nabla)u_2(t_1)\|_{DL_p(R^3)} \\
&\quad + \|(u_2(t_1) \bullet \nabla)u_2(t_1) - (u_2(t_2) \bullet \nabla)u_2(t_2)\|_{DL_p(R^3)} \\
&\leq M\left(2L + C_1|t_1 - t_2|^{\beta_1}\right)\|u_1(t_1) - u_2(t_1)\|_{DL_p(R^3)_{1/2}} \\
&\quad + 2\mathrm{MLC}|t_1 - t_2|^{\beta} \\
&\leq M(2L + C_1|t_1 - t_2|^{\beta_1})(\|u_1(t_1) - u_2(t_2)\|_{DL_p(R^3)_{\frac{1}{2}}} \\
&\quad + \|u_2(t_2) - u_2(t_1)\|_{DL_p(R^3)_{1/2}}) + 2\mathrm{MLC}|t_1 - t_2|^{\beta} \\
&\leq M(2L + 2^{\beta_1})\|u_1(t_1) - u_2(t_2)\|_{DL_p(R^3)_{1/2}} \\
&\quad + \mathrm{M}\ (2\mathrm{L}+2^{\beta_1})C_2|t_1 - t_2|^{\beta_2} + 2\mathrm{MLC}|t_1 - t_2|^{\beta} \\
&\leq L_1\|u_1(t_1) - u_2(t_2)\|_{DL_p(R^3)_{1/2}} + L_2|t_1 - t_2|^{\beta_3} \\
&\leq L_3(|t_1 - t_2|^{\beta_2} + \|u_1(t_1) - u_2(t_2)\|_{DL_p(R)_{1/2}}).
\end{aligned}$$

Hence, $F(t, u(t))$ satisfies the assumption $(F)$. Therefore, from lemma 4, for all initial data $(t_0, u_0) \in U$ the initial value problem (32) has a global solution

$$u(t) = \begin{cases} (t+1)(-\Delta)^{\frac{1}{2}}y_1(t), t \in [t_0, t_1] \\ (t+1)(-\Delta)^{\frac{1}{2}}y_2(t), t \in [t_1, \infty) \end{cases} \in C([t_0, \infty): DL_p(R^3)) \cap C^1((t_0, \infty): DL_p(R^3)),$$

i.e.,

$$u(t,x) = \begin{cases} (t+1)(-\Delta)^{-\frac{1}{2}}y_1(t) = T(t)u_0 + \int_{t_0}^{t} T(t-s)F(s, (-\Delta)^{-\frac{1}{2}}y_0(s))ds, t \in [t_0, t_1] \\ (t+1)(-\Delta)^{-\frac{1}{2}}y_2(t) = T(t)u_0 + \int_{t_0}^{t_1} T(t-s)F(s, (-\Delta)^{-\frac{1}{2}}y_0(s))ds, t \geq t_1\,. \end{cases} \qquad (37)$$

Using induction similar to Theorem 3.9 in [6], or, Theorem 5.1 in [21], we can prove that the solution $u(t,x) \in C^{\infty}([t_0, \infty) \times R^3)^3$.(see appendix) Substituting $u(t,x)$ into (1) we obtain solution $p(t,x)$. We also have $p(t,x) \in C^{\infty}([t_0, \infty) \times R^3)$. It follows from formula (2) and $u \in DL_p(R^3)$ that solution $u(t,x)$ is divergence-free. We also observe that $div\, u = 0$ from Theorem 1.6(ii) in [15]. Therefore, $u(t,x)$ and $p(t,x)$ are smooth global solutions of (1).

(Step 4) From the proof of Lemma 4, we see that the fixed point $y_0(t)$ in formula (21) belongs to $Y$, and the solution $u(t) = (-\Delta)^{-\frac{1}{2}} y_0(t)$. Therefore, $\|y_0(t)\|_{DL_p(R^3)} \leq M$ is bounded for $t \in [t_0, \infty)$ and some finite number $M$. Since $(-\Delta)^{-\frac{1}{2}}$ is bounded, then

$$\|u(t,x)\|_{DL_p(R^3)} = \left\|(-\Delta)^{-\frac{1}{2}} y_0(t)\right\|_{DL_p(R^3)} \leq \left\|(-\Delta)^{-\frac{1}{2}}\right\|_{B(DL_p(R^3))} M$$

for all $t \in [t_0, \infty)$and some constant $M$. Let $p = 2$,

$$\|u(t,x)\|_{DL_2(R^3)} = (\int_{R^3} |u(t,x)|^2 dx)^{1/2}$$

is bounded for all $t \in [t_0, \infty)$. Therefore, we have bounded energy.

To summarize the above, if the initial value $u_0$ is smooth and divergence free vector field, then the Navier-Stokes equation has a unique global smooth solution $p(t,x), u(t,x)$ on $[0,\infty) \times R^3$ with bounded energy.

Necessity. Suppose that the Navier-Stokes equation (1) has a unique smooth global solution $u$. Then $u$ satisfies (1), and so $u_0 \in C_{0,\sigma}^{\infty}(R^3)$. ☐

It is worth noting that in the proof of the above theorem, we did not use the condition (*). In fact, all the functions we discussed belong to $DL_p(R^3) \subseteq L_p(R^3)$, while the functions in $L_p(R^3)$ are all bounded in $R^3$.

Now let us discuss the specific form of the solution of N-S equation. If the Navier-Stokes equation has a unique smooth global solution $u(t,x)$ with bounded energy, then $u = (u_1, u_2, u_3)$ is smooth, divergence free and its three components have continuous partial derivatives. If every component $u_i(x_1,x_2,x_3)(i = 1,2,3)$ only changes in the direction of the axis $x_i$, so the rate of change (directional derivative) of $u_i(x_1,x_2,x_3)$ (i=1,2,3) is the largest in the direction of the axis $x_i$. We also know theoretically that the directional derivative of$u_i(x_1,x_2,x_3)(i$=1,2,3) is the largest in the direction of grad $u_i$= $(\partial u_i/\partial x_1, \partial u_i/\partial x_2, \partial u_i/\partial x_3)$. Therefore,

$$\partial_j u_i = 0 \quad (i \neq j, i,j = 1,2,3).$$

This result is obviously evident from an intuitive perspective. Then

$$du_i = \partial_i u_i \quad (i = 1,2,3).$$

It can be obtained by integration that $u_i = u_i(x_i)$. Since$u$is divergence free, $\sum_{i=1}^{3} \partial_i\, u_i(x_i) = 0$. By differentiating, we can obtain

$$\partial^2 u_i(x_i)/(\partial x_i)^2 = 0 (i = 1,2,3).$$

It can be seen from the integration that each $\partial_i u_i(x_i)$ can only be a constant and cannot contain $x_i$. It can be seen from further integration that each $u_i(x_i)$ can only contain linear term of $x_i$, i.e., $u_i(x_i) = C_i x_i + D_i$ where $C_i, D_i \in R$ $(i = 1,2,3)$ are constant. Thus, $u = (C_1 x_1 + D_1, C_2 x_2 + D_2, C_3 x_3 + D_3)$ does not satisfy bounded energy. Therefore, $C_i$=0$(i = 1,2,3)$ and $u = (D_1, D_2, D_3)$ is a uniform linear motion, i.e., $u = u(t,x)$ is a uniform linear motion when time $t$ is fixed. Specially, the initial value $u_0$ is a uniform linear motion. In fluid dynamics, a uniform flow is one in which all velocity vectors are identical (in both direction and magnitude) at every point of the flow for any given instant of time. Hence, the solution of N-S equation is a uniform flow in this case. But in our case, the flow is linear motion.

**Definition 4.** If a uniform flow is linear motion at every point, we call this type of uniform flow linear uniform flow.

Linear uniform flow is a special case of uniform flow, and is also a special case of laminal flow. Is this linear uniform flow the same linear motion at different times? From formula (37), it can be seen that this uniform flow is not a constant function of $t$, it is only a constant function of $x$ when time $t$ is fixed. If $t_1 \neq t_2$, then $T(t_1) \neq T(t_2)$, $T(t_1)u_0 \neq T(t_2)u_0$, the second integral items in the upper and lower formulas of the formula (37) are also different, hence $u(t_1, x) \neq u(t_2, x)$. This uniform flow is not only different constant vector field at the same point in space at different times, but the rate of change at the same point in space at different times may also be different. In fact,

$$\Delta\ u=0, (u \cdot \nabla)u =0$$

when $u = (u_1, u_2, u_3)$ is a linear uniform flow. From the formula (1) we have

$$\begin{pmatrix} \frac{du_1}{dt_1} \\ \frac{du_2}{dt_1} \\ \frac{du_3}{dt_1} \end{pmatrix} = \begin{pmatrix} \frac{\partial p(t_1, x)}{\partial x_1} \\ \frac{\partial p(t_1, x)}{\partial x_2} \\ \frac{\partial p(t_{1,}x)}{\partial x_3} \end{pmatrix} = \begin{pmatrix} \frac{\partial p(t_2, x)}{\partial x_1} \\ \frac{\partial p(t_2, x)}{\partial x_2} \\ \frac{\partial p(t_2, x)}{\partial x_3} \end{pmatrix} = \begin{pmatrix} \frac{du_1}{dt_2} \\ \frac{du_2}{dt_2} \\ \frac{du_3}{dt_2} \end{pmatrix}$$

According to Newton's First Law of Physics (the law of inertia), in the absence of external forces, an object will maintain its state of motion. But the interactions between molecules within the fluid consume energy, hence the velocity and acceleration may be different at different times and the same point in space.

Based on the above, we have obtained

**Theorem 2.** Suppose that the given externally applied force $f(t, x)$=0, and every component $u_i(x_1, x_2, x_3)(i = 1,2,3)$ of the solution $u(t, x)$ only changes in the direction of the axis $x_i$. Then this solution can only be a linear uniform flow.

## IV. CONCLUSION

There is a fundamental theorem in chaos mathematics which states that the development of anything depends on its initial state. Theorem 1 is consistent with this fundamental theorem. Theorem 1 tells us that as long as the initial value is smooth and divergence free vector field, then N-S equation has a unique global smooth solution. This complies with Newton's first law (the law of inertia), since we assume the external force is zero. The value of this solution at other times should not differ much from the initial value, but there will still be differences due to varying pressures in different locations.

We answered Millennium Prize Problems "Existence and smoothness of N-S solutions on $R^3$". The problem "Existence and smoothness of Navier-Stokes solutions in $R^3/Z^3$" has similar results. When $u_0$ is not smooth, divergence free vector field, then Navier-Stokes solutions may explode at some time or more solutions, which may be what Terence Tao, Tristan Buckmaster (Princeton) and Vlad Vicol (NYU) predict. This also answers the problems "Breakdown of Navier-Stokes solutions on $R^3$" and "Breakdown of Navier-Stokes solutions on $R^3/Z^3$".

## Appendix

**Lemma 1.** (Lemma 3.1 in [6]) Let u∈D $(-\Delta)$ and $-\Delta u \in (W^{m,r}(R^3))^3$ $(1< r < \infty)$ for some integer m $\geq 0$, then u∈ $(W^{m+2,r}(R^3))^3$ and satisfies

$\|u\|_{m+2,r} \leq C_{mr}\|-\Delta u\|_{m,r}$,

with a constant $C_{m,r} > 0$ independent of u and -Δ.

Let $C^{\mu}([0,\infty[;X)$ denote the space of Hölder continuous functions on $[0,\infty[$ with exponent $\mu$ and with values in a Banach space X. Similarly let $C^{\mu}((0,\infty);X)$ denote the space of functions which are Hölder continuous on every subinterval $[\varepsilon,\infty[$ of $(0,\infty)$, with exponent $\mu$.

**Lemma 2.** (Lemma 3.2 in [6]) Let $f(x) \in C^{\mu}\left([0,\infty[;DL_p(R^3)\right)$, for some $0 < \mu < 1$. Then the function

$$\varphi(t) = \int_0^t e^{-(t-s)(-\Delta)} f(s)ds \in C^{\mu}((0,\infty);D(-\Delta)) \cap C^{1+\tau}((0,\infty);DL_p(R^3))$$

for every $\tau$ such that $0< \tau < \mu$.

Let $P_r$ be the continuous projection from $L_p$ $(R^3)$ to $DL_p(R^3)$.

**Lemma 3.** (Lemma 3.3 in [6]) (i) For $m > 3/r$, there exists $C_{m,r}$ such that

$$\|P_r(u \cdot \nabla)\tau\|_{m,r} \leq C_{m,r}\|u\|_{m,r}\|\tau\|_{m+1,r}$$

for every $u \in (W^{m,r}(R^3))^3$, $\tau \in (W^{m+1,r}(R^3))^3$.

(ii) When r >3, we have

$\left\|P_r\ (u \cdot \nabla)\ \tau\right\|_{0,\ r} \leq C_r\|u\|_{1,\ r}\|\tau\|_{1,\ r}$.

for u,$\tau \in (W^{1,r}(R^3))^3$.

We will say that u(t) has property $(P_m)$(m $\geq$1) if

$u^{(m)} \in C^{\mu}((0,\infty);D(-\Delta)^{\frac{1}{2}})$,

$u^{(j)} \in C^{\mu}\left((0,\ \infty);\ \left(W^{m+1-j,2}(R^3)\right)^n\right), 1 \leq j \leq m-1$,

$u \in C^{\mu}\left((0,\infty);(W^{m+1,2}(R^3))^n\right)$,

for all $\mu$, $0 < \mu < \frac{1}{2}$. Here $u^{(j)} = \frac{d^j u}{dt^j}$.

Lemma 3(ii) implies

$-(u \cdot \nabla)u \in C^{\mu}((0,\infty);DL_p(R^3))$ for all $\mu, 0 < \mu < \frac{1}{2}$ .

Let the solution in the above Theorem 1

$$u(t,x)=\begin{cases} T(t)u_0+\int_0^t T(t-s)F\left(s,(-\Delta)^{-\frac{1}{2}}y_0(s)\right)ds, t \in [0,t_1], \\ T(t)u_0+\int_0^{t_1} T(t-s)F(s,(-\Delta)^{-\frac{1}{2}}y_0(s))ds, t \geq t_1. \end{cases}$$

Lemma 2 and Lemma 3(i) now imply

**Lemma 4.** (lemma 3.6 in [6]) $u \in C^{\mu}((0,\infty);\ D(-\Delta))$and $u' = {}^{du}/_{dt} \in C^{\mu}((0,\infty);\ \ DL_p(R^3))$ for all $\mu, 0 < \mu < {}^1/_2$.

Moreover,

$F(s,(-\Delta)^{-\frac{1}{2}}y_0(s)) \in C^{\mu}((0,\infty);\left(W^{1,2}(R^3)\right)^3)$.

**Lemma 5.** (Lemma 3.7 in [6]) We have u' $= {}^{du}/_{dx} \in C^{\mu}((0,\infty); D(\ (-\Delta)^{\frac{1}{2}}))$ for all $\mu, 0 < \mu < {}^1/_2$ .

The proof is similar to Lemma 3.7 in [6].

Since D$((-\Delta)^{\frac{1}{2}}$ ) ) is in $(W^{1,r}(R^3))^3$, Lemma 1, Lemma 3, Lemma 4 and the identity $u=(-\Delta)^{-1}(-(u\cdot\nabla)u-u)$ show that

$u \in C^{\mu}((0,\infty);(W^{3,2}(R^3)))^3$ .

Lemma 5 and the above formula show that u(t) has property$P_1$.

**Lemma 6.** (Lemma 3.8 in [6]) $P_m$ implies $P_1$.

The proof is the same to Lemma 3.8 in [6].

Therefore, we can prove the following theorem in a similar way as Theorem 3.9 in [6] or as Theorem 5.1 in [21].

**Theorem.** The solution in the above Theorem 1 is smooth.

We can also prove directly that the solution given by the above Theorem 1 is smooth. If we consider equation (32) in integral form

$$u(t)=e^{tP\Delta}u_0+\int_0^t e^{(t-s)P\Delta}\,Fu(s)ds. \tag{38}$$

Then solution u(t) in the above Theorem is the solution of (38). The Theorem 3.4 in [6] mean that as long as the solution of (38) exists, this solution is smooth. So according the Theorem 3.4 in [6] the solution in the above Theorem 1 is smooth.

**STATEMENT**

The purpose of publishing this article is to throw out a brick to attract jade, seeking the collective advice of mathematicians worldwide to tackle century-old mathematical problems.

The author declare that they have no conflict of interest. We have no any funding.

The datasets generated during the current study are available in the arXiv repository.

https://doi.org/10.48550/arXiv.2001.11699, the identifier is 2001.11699.

Authors' contributions

Maoting Tong conceived the entire article, read many literatures and proved the main Theorem1 and Theorem 2. Daorong Ton proved the relevant valuation lemma 3 and important lemma 4, and connect the reference literature. Maoting Tong and Daorong Ton have both solicited revision opinions from several peers. The results of several authors are cited in this article, and thanks are hereby given.

Maoting Tong was born in 1996, holds an undergraduate degree from Xi'anJiaotong Liverpool University and a PhD in Applied Mathematics from the University of Liverpool in UK. She joined Nanjing Vocational University of Technology in 2023, where her research focuses on Navier-Stokes equations.

Daorong Ton was born in January 1940. He ted from Wannan University in 1962 and has long beengradua engaged in mathematics teaching and research. From 1977 to 1987, he taught at the University of Science and Technology of China. From 1987 to 1988, he was a visiting professor at Boowling Green State University in the United States. After returning to China, he taught at Hohai University. In May 1995, he was a first-class visiting professor at the Maine University in France. He had published 58 mathematics papers. He has won awards multiple times and has chaired international mathematics conferences twice. Member of the Chinese Mathematical Society, member of the American Mathematical Society, member of the American Mathematical Association, reviewer for the American Mathematical Reviews, and judge for the selection committee of the ABI.